\documentclass[12pt]{amsart}
\usepackage{amsmath,amssymb,amsfonts,amsthm,amsopn}

\newcommand{\stft}{short-time Fourier transform}

\newcommand{\tf}{time-frequency}

\newcommand{\modsp}{modulation space}

\newtheorem{tm}{Theorem}[section]
\newtheorem{lemma}[tm]{Lemma}

\newtheorem{theorem}{Theorem}[section]

\newtheorem{proposition}[theorem]{Proposition}

\newcommand{\beqa}{\begin{eqnarray*}}
\newcommand{\eeqa}{\end{eqnarray*}}

\newcommand{\field}[1]{\mathbb{#1}}
\newcommand{\bR}{\field{R}}        
\newcommand{\bZ}{\field{Z}}        
\def\la{\lambda}

 \def\cF{\mathcal{F}}              
 \def\cS{\mathcal{S}}

 \def\cC{\mathcal{C}}

\def\rd{\bR^d}

\def\R{\right)}

\def\<{\left<}
\def\>{\right>}

\def\mv1{M_v^1}

\def\mpq{M^{p,q}}

\def\o{\xi}

\def\R{\mathbb{R}}
\def\Ren{\mathbb{R}^d}
\def\Renn{\mathbb{R}^{2d}}

\def\sch{\mathcal{S}}

\def\Fur{\mathcal{F}}

\def\f{\varphi}

\def\Sn2{S_{2}(L^{2}(\Ren))}
\def\S1{S_{1}(L^{2}(\Ren))}
\def\sig00{\sigma_{0,0}}

\def\la{\langle}
\def\ra{\rangle}

\newcommand{\wpq}{W(L^p,L^q)}

\begin{document}

\title[]{Sharpness of some properties of Wiener amalgam and modulation spaces}

\author{Elena Cordero and Fabio Nicola}
\address{Department of Mathematics,  University of Torino,
Via Carlo Alberto 10, 10123
Torino, Italy}
\address{Dipartimento di Matematica, Politecnico di
Torino, Corso Duca degli
Abruzzi 24, 10129 Torino,
Italy}
\email{elena.cordero@unito.it}
\email{fabio.nicola@polito.it}



\date{}

\begin{abstract} We prove
sharp estimates for the
dilation operator
$f(x)\longmapsto f(\lambda
x)$, when acting on Wiener
amalgam spaces $W(L^p,L^q)$.
Scaling arguments are also used
to prove the sharpness of the
known convolution and
pointwise relations for
modulation spaces $M^{p,q}$,
as well as the optimality of
an estimate for the
Schr\"odinger propagator on
modulation spaces.
\end{abstract}

\maketitle

\section{Introduction}
Modulation and Wiener amalgam
spaces have been introduced
and used to measure the \tf\
concentration of functions
and tempered distributions in
the framework of
time-frequency analysis
\cite{feichtinger83,F1,feichtinger90,feichtinger-norbert,book,Heil03,triebel83}.
Recently, these spaces have
been employed to study
boundedness properties of
pseudodifferential operators
(see, e.g.,
\cite{CG02,sugitomita2,Toft04}),
Fourier Integral operators
(in particular, Fourier
multipliers)
\cite{benyi,fio1,cordero-nicola-rodino}
and wellposedness of
solutions to PDE's (see,
e.g.,
\cite{benyi2,benyi3,cordero2,cordero3,baoxiang3,baoxiang2,baoxiang}
and references therein).

In this paper we present
new dilation properties for
Wiener amalgam spaces and their optimality. Moreover,  we
prove the sharpness of the
known convolution and
pointwise estimates for
modulation spaces.

To recall  the definition of
these spaces, we first  introduce the translation and
modulation operators, defined   by
$
 T_xf(t)=f(t-x)\quad{\rm and}\quad M_{\o}f(t)= e^{2\pi i \o
  t}f(t) \,
$,  $t,x,\xi\in\rd$.
\par
\textit{Wiener amalgam
spaces}
\cite{F1,feichtinger90,Heil03}.
Let $g \in \cC_0^\infty(\rd )$ be a
test function. We will refer
to $g$ as a window function.
Let $B$ be either
 the Banach space
$L^p(\rd )$ or  $\cF L^p(\rd )$,
 $1\leq p\leq\infty$. For any given function $f$
which is locally in $B$ (i.e.
$g f\in B$, $\forall
g\in\cC_0^\infty(\rd )$), we set
$f_B(x)=\| fT_x g\|_B$. The
Wiener amalgam space
$W(B,L^q)(\rd )$ with local
component $B$ and global
component  $L^q(\rd )$, $1\leq
q\leq\infty$, is defined as
the space of all functions
$f$ locally in $B$ such that
$f_B\in L^q(\rd )$. Endowed with
the norm
$\|f\|_{W(B,L^q)}=\|f_B\|_{L^q}$,
$W(B,L^q)(\rd )$ is a Banach space.
Moreover, different choices
of $g\in \cC_0^\infty(\rd )$
generate the same space and
yield equivalent norms. In
fact, the space of admissible
windows for the Wiener
amalgam spaces $W(B,L^q)(\rd )$ can
be enlarged to the so-called
Feichtinger algebra $W(\Fur
L^1,L^1)(\rd )$. Recall that the
Schwartz class $\sch(\rd )$  is
dense in $W(\Fur
L^1,L^1)(\rd )$.\par

\textit{Modulation spaces}
\cite{F1, book}.
For a fixed non-zero  $g \in \cS (\rd )$  the \stft\ of $f \in
 \cS ' (\rd ) $ with respect to the window $g$ is given by
$   V_gf(x,\o)=\la f,M_\o T_x g\ra $, \,$x,\xi\in\rd$.\par
 Given a
non-zero window
$g\in\sch(\Ren)$, $1\leq
p,q\leq \infty$, the {\it
  modulation space} $M^{p,q}(\Ren)$ consists of all tempered
distributions $f\in\sch'(\Ren)$ such that $V_gf\in L^{p,q}(\Renn
)$ (mixed-norm spaces). The norm on $M^{p,q}$ is given by
$$
\|f\|_{M^{p,q}}:=\|V_gf\|_{L^{p,q}}=\left(\int_{\Ren}
  \left(\int_{\Ren}|V_gf(x,\o)|^p\,
    dx\right)^{q/p}d\o\right)^{1/q}  \,
    ,
$$
 with
obvious changes if $p=\infty$
or $q=\infty$. If $p=q$, we
write $M^p$ instead of
$M^{p,p}$. The space
$M^{p,q}(\Ren )$ is a Banach
space, whose definition is
independent of the choice of
the window $g$. Different
non-zero windows $g\in M^1$
yield equivalent norms on
$\mpq$. This property will be
crucial in the sequel,
because we will choose a
suitable window $g$ in
estimates of the
$M^{p,q}$-norm.  Within the
class of \modsp s, one finds
several standard function
spaces, for instance $M^2 =
L^2$, $M^1=W(\cF
L^1,L^1)$ and, using weighted
versions, one also finds
certain Sobolev spaces and
Shubin-Sobolev spaces
(\cite{CG02,book}). The
relationship between
modulation and Wiener amalgam
spaces is expressed by  the
following
result:\par\vskip0.1truecm
\textsl{The Fourier transform
establishes an isomorphism
$\Fur: M^{p,q}\to W(\Fur
L^p,L^q)$.}\vskip0.1truecm
Consequently, convolution
properties of modulation
spaces can be translated into
pointwise multiplication
properties of Wiener amalgam
spaces.\par\vskip0.1truecm
Let us now turn to the topic
of the present paper. The
importance of the dilation
operator $$f(x)\longmapsto
f(\lambda x),\quad
\lambda>0,$$ in classical
analysis  is well-known. For
example, in most estimates
arising in classical Harmonic
Analysis (e.g. the H\"older,
Young and Hausdorff-Young
inequalities) as well as in
Partial Differential
Equations (e.g. Sobolev
embeddings, Strichartz
estimates) scaling arguments
yield the constraints that
the Lebesque exponents must
satisfy for the corresponding
inequalities to hold.
\par When dealing with
modulation or Wiener amalgam
spaces, the situation becomes
more subtle. In fact, the
corresponding norms are not
`homogeneous' with respect to
the scaling. Basically, this
is due to the fact that, for
example, in $W(L^p,L^q)$ the
two spaces $L^p$, $L^q$
display different scaling if
$p\not=q$. Obtaining sharp
estimates (in terms of
$\lambda$) for  the dilation
operator  norm, when acting
on such spaces, is therefore
a non trivial problem. This
study was widely carried out
in \cite{sugimototomita} (see
also \cite{cordero2,Toft04})
in the case of modulation
spaces $M^{p,q}$. The
estimates obtained in
\cite{sugimototomita} turned
out to be a fundamental tool
for embedding problems of
modulation spaces into Besov
spaces  (see also
\cite{baoxiang2}), and for
boundedness of
pseudodifferential operators
of type $(\rho,\delta)$ on
modulation spaces
(\cite{sugitomita2}). Let us
highlight that this type of
arguments allowed us to prove
the sharpness of some
Strichartz estimates in the
Wiener amalgam spaces $W(\Fur
L^p,L^q)$ in \cite{cordero3}.
Finally, they were also used
in
\cite{cordero-nicola-rodino}
to prove sharp boundedness
properties of H\"ormader's
type Fourier integral
operators on $\Fur L^p$ and
modulation spaces.

Let us point out
that an investigation of the
dilation operator on
$W(C,L^1)$ ($C$ being the space of continuous functions) had already
appeared in
\cite{feichtinger-norbert}.\par
 The first result of this note (Proposition \ref{dil} below) provides {\it
 sharp}
upper and lower bounds for
the operator norm of the
dilation operator on the
Wiener amalgam spaces
$W(L^p,L^q)$. \par
Differently from what  one
could ingenuously expect, it
does not happen that the
exponent $p$ alone has its
influence when
$\lambda\to+\infty$, whereas
the exponent $q$ when
$\lambda\to0$.

Then, as for the classical function spaces, scaling arguments
are employed to prove the
sharpness of the known
 convolution, inclusion and  pointwise
moltiplication relations for modulation
spaces. This is precisely the
topic studied in Section $3$. To chase this goal,
we do not use the bounds obtained in
\cite{sugimototomita}, that  would
give weaker constraints than
the optimal ones. Instead, the sharp results are the issues of
 explicit computations involving dilation properties of
Gaussian functions.
\par Finally, we observe that these techniques can be applied to prove
 the sharpness of  estimates arising in PDEs. As an example, in Section 4 we
prove the optimality of an
estimate for the
Schr\"odinger propagator,
recently obtained in
\cite{baoxiang2} (see also
\cite{benyi}).\par

\section{Dilation properties
of Wiener amalgam
spaces}\label{dilWiener} In
this section we study the
dilation properties of Wiener
amalgam spaces $W(L^p,L^q)$, $1\leq p,q\leq\infty$.
First, recall the following complex
interpolation result \cite{feichtinger83}.
\begin{lemma}\label{interpolazione}
  Let $B_0,B_1$, be local components
  of Wiener amalgam spaces,
   as in the Introduction.
   Then, for $1\leq
   q_0,q_1\leq\infty$, with
   $q_0<\infty$ or
   $q_1<\infty$,
and $0<\theta<1$, we
  have
\[
  [W(B_0,L^{q_0}),W(B_1,L^{q_1})]_{[\theta]}
  =W\left([B_0,B_1]_{[\theta]},L^{q_\theta}\right),
  \]
with
$q_\theta=(1-\theta)/{q_0}+\theta/{q_1}$.
  \end{lemma}

 For $\lambda>0$, we set
$f_\lambda(x)=f(\lambda x)$.
\begin{proposition}\label{dil}
For $1\leq p,
q\leq\infty$,
\begin{equation}\label{dil1}
\|f_\lambda\|_{W(L^p,L^q)}\lesssim
\lambda
^{-d\max\left\{\frac{1}{p},\frac{1}{q}\right\}}\|f\|_{\wpq},\quad
\forall\, 0<\lambda\leq1,
\end{equation}
and
\begin{equation}\label{dil2}
\|f_\lambda\|_{W(L^p,L^q)}\lesssim
\lambda
^{-d\min\left\{\frac{1}{p},\frac{1}{q}\right\}}\|f\|_{\wpq},\quad
\forall \lambda\geq1.
\end{equation}
Also, we have
\begin{equation}\label{dil1m}
\|f_\lambda\|_{W(L^p,L^q)}\gtrsim
\lambda ^{-d\min
\left\{\frac{1}{p},\frac{1}{q}\right\}}\|f\|_{\wpq},\quad
\forall\, 0<\lambda\leq1,
\end{equation}
and
\begin{equation}\label{dil2m}
\|f_\lambda\|_{W(L^p,L^q)}\gtrsim
\lambda ^{-d\max
\left\{\frac{1}{p},\frac{1}{q}\right\}}\|f\|_{\wpq},\quad
\forall \lambda\geq1.
\end{equation}
\end{proposition}
We first prove the following
weaker estimates.
\begin{lemma}
For $1\leq p,
q\leq\infty$,
\begin{equation}\label{dil1bis}
\|f_\lambda\|_{W(L^p,L^q)}\lesssim
\lambda
^{-d\left(\frac{1}{p}+\frac{1}{q}\right)}
\|f\|_{\wpq},\quad \forall\,
0<\lambda\leq1,
\end{equation}
and
\begin{equation}\label{dil2bis}
\|f_\lambda\|_{W(L^p,L^q)}\lesssim
\lambda
^{d\left(1-\frac{1}{p}-\frac{1}{q}\right)}
\|f\|_{\wpq},\quad
\forall\lambda\geq1.
\end{equation}
\end{lemma}
\begin{proof}
To compute the Wiener norm, we choose the window function $g=\chi_{B(0,1)}$, the characteristic
function of the ball $B(0,1)$. Then,
\begin{eqnarray*}
\|f_\lambda\|_{\wpq}&\asymp&\|\|f(\lambda
t)
g(t-x)\|_{L^p_t}\|_{L^q_x}\\
&=&\lambda^{-\frac{d}{p}}\|\|f(t)g_{1/\lambda}(t-\lambda
x)\|_{L^p_t}\|_{L^q_x}\\
&=&\lambda
^{-d\left(\frac{1}{p}+\frac{1}{q}\right)}\|\|f(t)g_{1/\lambda}(t-
x)\|_{L^p_t}\|_{L^q_x}.
\end{eqnarray*}
If $0<\lambda\leq 1$, the window function $g$ fulfills
$g_{1/\lambda}(y)\leq g(y)$,
and \eqref{dil1bis} follows.\par
\noindent To prove \eqref{dil2bis}, we
argue by duality. Indeed, if
$\lambda\geq1$, relation
 \eqref{dil1bis}, applied
to the pair $(p',q')$, yields
\begin{eqnarray*}
\|f_\lambda\|_{\wpq}&=&\sup_{\|g\|_
{W(L^{p'},L^{q'})}=1}|\langle
f_\lambda,g\rangle|\\
&=&
\sup_{\|g\|_{W(L^{p'},L^{q'})}=1}
\lambda^{-d}|\langle
f,g_{1/\lambda}\rangle|\\
&\leq&
\lambda^{-d}\sup_{\|g\|_{W(L^{p'},L^{q'})}=1}
\|f\|_{\wpq}\|g_{1/\lambda}
\|_{W(L^{p'},L^{q'})}\\
 &\lesssim&
\lambda^{-d}\lambda^{d\left(\frac{1}{p'}
+\frac{1}{q'}\right)}\|f\|_{\wpq}.
\end{eqnarray*}
\end{proof}
\begin{proof}[Proof of
Proposition \ref{dil}] We
first prove \eqref{dil1} and
\eqref{dil2} when $p=\infty$.
We see at once that
\eqref{dil1} coincides with
\eqref{dil1bis} when
$p=\infty$. On the other
hand, \eqref{dil2} for
$p=\infty$ follows by complex
interpolation (Lemma
\ref{interpolazione}) from
\eqref{dil2bis} with
$(p,q)=(\infty,1)$, i.e.,
\[
\|f_\lambda\|_{W(L^\infty,L^1)}\lesssim
\|f\|_{W(L^\infty,L^1)}
\]
and the trivial estimate
\[
\|f_\lambda\|_{W(L^\infty,L^\infty)}\asymp\|f_\lambda\|_{L^\infty}=\|f\|_{L^\infty}\asymp\|f\|_{W(L^\infty,L^\infty)}.
\]
Since the estimates
\eqref{dil1} and \eqref{dil2}
also hold for $p=q$ (because
$W(L^p,L^p)=L^p$ with
equivalent norms), by
interpolation with the case
$p=\infty$, $1\leq q\leq
\infty$, we see that they
hold for any pair $(p,q)$,
with $1\leq q\leq
p\leq\infty$. When $p<q$
they follow by duality  arguments as in the
proof of \eqref{dil2bis}.\par
Finally, \eqref{dil1m} and
\eqref{dil2m} follow at once
from \eqref{dil2} and
\eqref{dil1}, respectively, applied to the function
$f_{1/\lambda}$.\par
\end{proof}\par
We now show that the
result above is sharp.
\begin{proposition}
{\rm (Sharpness of
\eqref{dil1} and \eqref{dil2})}. \\
(i) Suppose that, for some
$\alpha\in\R$,
\begin{equation}\label{dil1bis2}
\|f_\lambda\|_{W(L^p,L^q)}\lesssim
\lambda
^{\alpha}\|f\|_{\wpq},\quad
\forall\, 0<\lambda\leq1.
\end{equation}
Then
\begin{equation}
\alpha\leq
-d\max\left\{\frac{1}{p},\frac{1}{q}\right\}.
\end{equation}
(ii) Suppose that, for some
$\alpha\in\R$,
\begin{equation}\label{dil2bis2}
\|f_\lambda\|_{W(L^p,L^q)}\lesssim
\lambda
^{\alpha}\|f\|_{\wpq},\quad
\forall \lambda\geq1.
\end{equation}
Then \begin{equation}
\alpha\geq
-d\min\left\{\frac{1}{p},\frac{1}{q}\right\}.
\end{equation}
\end{proposition}
This also shows the sharpness
of the estimates
\eqref{dil1m} and
\eqref{dil2m}, since they
are equivalent to
\eqref{dil2} and \eqref{dil1},
respectively.
\begin{proof}
\textit{(i)} First, consider  the case
$p\geq q$. We have
$W(L^p,L^q)\hookrightarrow
W(L^q,L^q)=L^q$.
Hence
\[
\lambda^{-\frac{d}{q}}\|f\|_{L^q}=\|f_\lambda\|_{L^q}\lesssim
\|f_\lambda\|_{W(L^p,L^q)}.
\]
Combining this estimate with
\eqref{dil1bis2} and letting
$\lambda\to0^+$, we obtain
$\alpha\leq -d/q$.\par
Assume now  $p<
q$. It suffices to verify
that, for every $\epsilon>0$,
there exists $f\in
W(L^p,L^q)$ such that
\begin{equation}\label{ese}
\|f_\lambda\|_{W(L^p,L^q)}\geq
C\lambda^{-\frac{d}{p}+\epsilon}.
\end{equation}
We study the case of
dimension $d=1$. The general
case follows by tensor
products of functions of one
variable. To this end, we
choose
\[
f(t)=\begin{cases}|t|^{-\frac{1}{p}+\epsilon}\
&{\rm for}\ |t|\leq1\\
0& {\rm for}\ |t|>1.
\end{cases}
\]
Observe that  $f\in
W(L^p,L^1)\hookrightarrow
W(L^p,L^q)$, for every $1\leq
q\leq\infty$, and
\begin{equation}\label{due}
f(\lambda
t)=\lambda^{-\frac{1}{p}+\epsilon}f(t),\quad{\rm
for}\
|t|\leq\frac{1}{\lambda}.
\end{equation}
Now,  take $g=\chi_{B(0,1)}$ as
window function. Of course,
\begin{align*}
\|f_\lambda\|_{W(L^p,L^q)}&=\left(\int
\|f_\lambda T_y g\|_{L^p}^q
dy\right)^{1/q}\\
&\geq\left(\int_{B(0,1)}
\|f_\lambda T_y g\|_{L^p}^q
dy\right)^{1/q}.
\end{align*}
By using \eqref{due} and the choice  $g=\chi_{B(0,1)}$, for $\lambda\leq 1/2$
the last expression is estimated from below by
\[
\geq
\lambda^{-1/p+\epsilon}\left(\int_{B(0,1)}
\|f T_y g\|_{L^p}^q
dy\right)^{1/q},
\]
that is \eqref{ese}.\par $(ii)$
Again, we first consider
the case $p\geq q$, namely
$q'\geq p'$. Then
$L^{p'}=W(L^{p'},L^{p'})\hookrightarrow
W(L^{p'},L^{q'})$. Hence,
\[
\|f_\lambda\|_{W(L^p,L^q)}=\sup_{\|g\|_{W(L^{p'},L^{q'})}=1}|\langle
f_\lambda, g\rangle|\gtrsim
\sup_{\|g\|_{L^{p'}=1}}|\langle
f_\lambda,
g\rangle|=\|f_\lambda\|_{L^p}=\lambda^{-\frac{d}{p}}\|f\|_{L^p}.
\]
Combining this estimate with
\eqref{dil2bis2} and letting
$\lambda\to+\infty$, we obtain
$\alpha\geq-\frac{d}{q}$.\par
Suppose now $p<q$. As before
it suffices to prove, in
dimension $d=1$, that for
every $\epsilon>0$ there
exists a function $f\in
W(L^p,L^q)$ such that
\[
\|f_\lambda\|_{W(L^p,L^q)}\geq
C\lambda^{-\frac{1}{q}-\epsilon}.
\]
Therefore, choose
\[
f(t)=\begin{cases}|t|^{-\frac{1}{q}-\epsilon}\
&{\rm for}\ |t|\geq1\\
0& {\rm for}\ |t|<1.
\end{cases}
\]
Then $f\in
W(L^\infty,L^q)\hookrightarrow
W(L^p,L^q)$, for every $1\leq
p\leq\infty$, and
\begin{equation}\label{tre}
f(\lambda
t)=\lambda^{-\frac{1}{q}-\epsilon}f(t),\quad{\rm
for}\
|t|\geq\frac{1}{\lambda}.
\end{equation}
Again, choose $g=\chi_{B(0,1)}$
as  window function. We have
\begin{align*}
\|f_\lambda\|_{W(L^p,L^q)}\geq\left(\int_{B(0,2)}
\|f_\lambda T_y g\|_{L^p}^q
dy\right)^{1/q}.
\end{align*}
By using \eqref{tre} and
 the choice $g=\chi_{B(0,1)}$, for
$\lambda\geq 1$
the last expression  is
\[
\geq
\lambda^{-1/q-\epsilon}\left(\int_{B(0,2)}
\|f T_y g\|_{L^p}^q
dy\right)^{1/q},
\]
which concludes the proof of $(ii)$.
\end{proof}
\par

\section{Convolution,  inclusion and multiplication relations for modulation spaces}
In this section we study the
optimality of the convolution,
inclusion and pointwise multiplication relations for
modulation spaces. We need some preliminary results.

If one chooses the Gaussian $e^{-\pi |x|^2}$ as window
 function to compute   Wiener amalgam  norms,
 then an easy computation (see e.g.
  \cite[Lemma 5.3]{cordero3}) yields the  result below.
\begin{lemma}\label{lemm2}
For $a,b\in\R$, $a>0$, set
$\mathcal{G}_{(a+ib)}(x)=(a+ib)^{-d/2}e^{-\frac{\pi|x|^2}{a+ib}}$.
Then, for every $1\leq
p,q\leq\infty$,
\begin{equation}\label{lem2}
\|\mathcal{G}_{(a+ib)}\|_{W(\Fur
L^p,L^q)}=\frac{\left((a+1)^2+b^2\right)^{
\frac{d}{2}\left(\frac{1}{p}-\frac{1}{2}
\right)}}{p^\frac{d}{2p}(aq)^\frac{d}{2q}
\left(a(a+1)+b^2\right)^{\frac{d}{2}\left(\frac{1}{p}
-\frac{1}{q}\right)}}.
\end{equation}\end{lemma}
For tempered distributions compactly supported either in time or in frequency, the $M^{p,q}$-norm is equivalent  to
 the $\cF L^q$-norm or $L^p$-norm, respectively. This result is well-known (\cite{fe89-1,feichtinger90,kasso07}). For the sake of
completeness we provide an
outline of the proof.
\begin{lemma}\label{lloc} Let $1\leq p,q\leq
\infty.$\\
(i) For every $u\in
\mathcal{S}'(\rd)$, supported
in a compact set $K\subset
\rd$, we have $u\in
M^{p,q}\Leftrightarrow u\in
\Fur L^q$, and
\begin{equation}\label{loc}
C_K^{-1} \|u\|_{M^{p,q}}\leq
\|u\|_{\cF L^q}\leq C_K
\|u\|_{M^{p,q}},
\end{equation}
where $C_K>0$ depends only on
 $K$.\\
 (ii) For every $u\in \mathcal{S}'(\rd)$,
whose Fourier transform is
supported in a compact set
$K\subset \rd$, we have $u\in
M^{p,q}\Leftrightarrow u\in
L^p$, and
\begin{equation}\label{loc2}
C_K^{-1} \|u\|_{M^{p,q}}\leq
\|u\|_{ L^p}\leq C_K
\|u\|_{M^{p,q}},
\end{equation}
where $C_K>0$ depends only on
 $K$.\\
\end{lemma}
\begin{proof}
$(i)$ It is detailed in  \cite[Lemma 1]{kasso07}.\\
\noindent $(ii)$ It is well-known (see e.g.
\cite{triebel83}) that
\[
\|u\|_{M^{p,q}}\asymp\left(\sum_{k\in\bZ^d}\|\nu(D-k)
u\|_{L^p}^q\right)^{1/q},
\]
where $\nu$ is a test
function satisfying
$\sum_{k\in\bZ^d}\nu(\xi-k)\equiv1$.
Now, if $\hat{u}$ has compact
support, the above sum is
finite and one deduces at
once the first estimate in
\eqref{loc2}, since the
multipliers $\nu(D-k)$ are
(uniformly) bounded on $L^p$.
To obtain the second estimate
in \eqref{loc2}, we write
$u=\sum_{k\in \bZ^d}\nu(D-k)
u$, then apply the
triangle inequality and
the finiteness of the sum over
$k$ again.
\end{proof}

Now, we turn our attention to the sharpness of the convolution properties for \modsp s.

\begin{proposition}\label{mconvmp}
 Let $1\leq p,q,p_1,p_2,q_1,q_2\leq\infty$. Then
\begin{equation}\label{mconvm}
\|f\ast g\|_{M^{p,q}}\lesssim \|f\|_{ M^{p_1,q_1}}\|g\|_{M^{p_2,q_2}}
\end{equation}
if and only if the following
indices' relations hold true:
\begin{equation}\label{indyou}
\frac1p+1\leq\frac1{p_1}+\frac1{p_2},
\end{equation}  and
\begin{equation}\label{indho}
\frac1{q}\leq \frac1{q_1}+\frac1{q_2}.
\end{equation}
\end{proposition}
\begin{proof}\emph{Sufficiency.} The inclusion relations
\eqref{mconvm},  were proved
in \cite{CG02,Toft04}. There
the indices' relations
\eqref{indyou} and
\eqref{indho} were shown with
the equalities. The
inequalities follow by the
inclusion relations
$M^{p_1,q_1}\hookrightarrow
M^{p_2,q_2}$ for $p_1\leq
p_2$ and $q_1\leq q_2$
(\cite{F1,book}).\par
\noindent \emph{Necessity.}
We consider the family of
Gaussians
 $\f_\lambda(x) :=e^{-\pi\lambda|x|^2}$, for
$\lambda>0$. Obviously,
$\f_\lambda\in\cS(\Ren)\subset
M^{p,q}(\Ren)$, for every
$1\leq p,q\leq\infty$. Since
$\|f\|_{M^{p,q}}\asymp
\|\hat{f} \|_{W(\cF
L^p,L^q)}$ and
$\hat{\f}_\lambda=\lambda^{-d/2}\f_{1/\lambda}$,
Lemma \ref{lemm2} yields:
\begin{equation}\label{mpqGauss}
   \|\f_\lambda\|_{M^{p,q}}\asymp \lambda^{-d/2}\|
   \f_{1/\lambda}\|_{W(\cF L^p,L^q)}\asymp
   \|\mathcal{G}_{(\lambda)}\|_{W(\cF L^p,L^q)}\asymp
   \frac{(\lambda+1)^{d(1/p-1/2)}}{\lambda^{d/(2q)}(\lambda^2+\lambda)^{(1/p-1/q)d/2}}.
\end{equation}
A straightforward calculation
shows that $(\f_\lambda\ast
\f_\lambda)(x)=(2\lambda)^{-\frac
d2}\f_{\lambda/2}(x)$. Hence,
using \eqref{mpqGauss}, we
obtain
\begin{equation}\label{estres}
\|\f_\lambda\ast
\f_\lambda\|_{M^{p,q}}\asymp\lambda^{-(1+1/p)d/2},\quad\mbox{for}\,\,\lambda\rightarrow
0^+.
\end{equation}
Using \eqref{mpqGauss} again, we also obtain
\begin{equation}\label{bo}
\|\f_\lambda\|_{M^{p_i,q_i}}\asymp\lambda^{-\frac
d{2p_i}},\quad
i=1,2,\qquad\mbox{for}\,\,\lambda\rightarrow
0^+.
\end{equation}
Substituting in \eqref{mconvm}, we obtain \eqref{indyou}. The relation \eqref{indho} can be obtained
similarly. Indeed, the estimate
\eqref{mpqGauss} gives, for $\lambda\rightarrow+\infty$,
$$\|\f_\lambda\ast
\f_\lambda\|_{M^{p,q}}\asymp\lambda^{-d(1-\frac
1{2 q})},\quad
\|\f_\lambda\|_{M^{p_i,q_i}}\asymp\lambda^{-\frac
d2(1-\frac 1{q_i})},\quad
i=1,2,$$ and, using
\eqref{mconvm} again, the
relation \eqref{indho}
follows.\par An alternative
proof of the necessary
conditions \eqref{indyou} and
\eqref{indho} is provided by
Lemma \ref{lloc}.
 Precisely, to
prove \eqref{indho}, consider
two compactly supported
smooth functions $f,g$ and their scaling
$f_\lambda(x)=f(\lambda x)$, $g_\lambda(x)=g(\lambda
x)$, with $\lambda\geq1$. Since  $\lambda\geq1$, $f_\lambda$ and
$g_\lambda$ (and therefore
$f_\lambda\ast g_\lambda$)
are all supported in a compact
subset $K$, independent of
$\lambda$. By Lemma \ref{lloc},  $(i)$,
 the bilinear estimate \eqref{mconvm} for $f_\lambda$ and
$g_\lambda$  becomes
$$\|f_\lambda\ast g_\lambda\|_{\cF L^q}\lesssim \|f_\lambda\|_{\cF L^{q_1}}\|g_{\lambda}\|_{\cF L^{q_2}}.
$$
Using $f_\lambda\ast
g_\lambda=\lambda^{-d}(f\ast
g)_\lambda$, the dilation
property for $\cF L^q$
spaces:
$\|h(\lambda\cdot)\|_{\cF
L^q}=\lambda^{-\frac
d{q'}}\|h\|_{\cF L^q}$,
 and letting
$\lambda\to+\infty$, we obtain
\eqref{indho}.\par
In order to prove
\eqref{indyou}, one argues
similarly.  Here  the functions  $f,g$ have  Fourier transforms $\hat{f}, \hat{g}$
 compactly supported and  the scale $\lambda$ satisfies $0<\lambda\leq1$. By Lemma
\ref{lloc}, $(ii)$, the estimate \eqref{mconvm}  becomes
$$\|f_\lambda\ast g_\lambda\|_{ L^p}\lesssim \|f_\lambda\|_{ L^{p_1}}\|g_{\lambda}\|_{ L^{p_2}}.
$$
Using $f_\lambda\ast g_\lambda=\lambda^{-d}(f\ast g)_\lambda$, the dilation  property $\|h(\lambda\cdot)\|_{ L^p}=\lambda^{-d/p}\|h\|_{ L^p}$, and letting
$\lambda\to0^+$, we prove \eqref{indyou}.
\end{proof}
\par
\vskip0.3truecm
The family of Gaussians
$\f_\lambda$ provides an
alternative proof for the
sharpness of the inclusion
relation for \modsp s, already
 obtained by the
inclusion relations for the
sequence spaces $\ell^{p,q}$,
via the norm equivalence
$\|f\|_{M^{p,q}}\asymp \|\la
f,T_mM_n
g\ra\|_{\ell^{p,q}}$, with
$\{ T_mM_n g\}$ being a
 Gabor frame (see, e.g., \cite[Theorem 13.6.1]{book}).
\begin{proposition} Let
$1\leq
p_1,p_2,q_1,q_2\leq\infty$.
Then,
\begin{equation}\label{incl1}\|f\|_{M^{p_2,q_2}}\lesssim \|f\|_{M^{p_1,q_1}}
\end{equation}
if and only if the following indices' relation holds
\begin{equation}\label{incl2}
p_1\leq p_2 \quad \mbox{and}\quad \,q_1\leq q_2.
\end{equation}
\end{proposition}
\begin{proof} We show
the necessity of \eqref{incl2}. Let
$\f_\lambda(x)
=e^{-\pi\lambda|x|^2}$,
$\lambda>0$. From the
proof of Proposition
\ref{mconvmp},
\begin{equation}\label{ABCD}
\|\f_\lambda\|_{M^{p_i,q_i}}\asymp\lambda^{-\frac
d{2p_i}},\quad\mbox{for}\,\lambda\rightarrow
0^+\quad\mbox{and}\quad\|\f_\lambda\|_{M^{p_i,q_i}}\asymp\lambda^{-\frac
d{2}(1-\frac1{q_i})},\quad\mbox{for}\,\lambda\rightarrow
+\infty. \end{equation}
Hence, for \eqref{incl1} to
be satisfied it must be
$$\lambda^{-\frac
d{2p_2}}\lesssim
\lambda^{-\frac
d{2p_1}},\quad\mbox{for}\,\lambda\rightarrow
0^+\quad\mbox{and}\quad\lambda^{-\frac
d2(1-\frac1{q_2})}\lesssim
\lambda^{-\frac
d2(1-\frac1{q_1})},\quad\mbox{for}\,\lambda\rightarrow+\infty,
$$ that give the indices' relations in \eqref{incl2}.
\end{proof}

\vskip0.3truecm
In what follows we study  the
pointwise multiplication operator in modulation spaces
(which is equivalent to studing  the convolution operator
for the Wiener amalgam spaces $W(\Fur
L^p,L^q)$).
\begin{proposition}\label{mconvmp2}
 Let $1\leq  p,q,p_1,p_2,q_1,q_2\leq\infty$. Then
\begin{equation}\label{mconvm2}
\|f g\|_{M^{p,q}}\lesssim
\|f\|_{M^{p_1,q_1}}\|g\|_{M^{p_2,q_2}}
\end{equation}
if and only if the following
indices' relations hold true:
\begin{equation}\label{indyou2}
\frac1p\leq\frac1{p_1}+\frac1{p_2},
\end{equation}  and
\begin{equation}\label{indho2}
\frac1{q}+1\leq
\frac1{q_1}+\frac1{q_2}.
\end{equation}
\end{proposition}
\begin{proof} The sufficiency
can be found in \cite{F1}
(see also \cite{baoxiang}).
For the necessity of the
conditions \eqref{indyou2}
and \eqref{indho2} we test
the estimate \eqref{mconvm2}
on the Gaussians
$f(x)=g(x)=\f_\lambda(x)=e^{-\lambda\pi
|x|^2}$. We observe that
$\f_\lambda
\f_\lambda=\f_{2\lambda}$.
Hence by applying
\eqref{ABCD} and substituting
in \eqref{mconvm2}, relation
\eqref{indho2} follows by
letting $\lambda\to 0^+$,
whereas \eqref{indyou2}
follows by letting
$\lambda\to+\infty$.
\end{proof}
\section{An estimate for the
Schr\"odinger propagator}
Consider the Fourier multiplier $e^{it\Delta}$, with
symbol $e^{-it|2\pi \xi|^2}$,
i.e.,
 \[
 \left(e^{it\Delta} u_0\right)(x)=\frac{1}{(4\pi i
t)^{d/2}} \int
e^{i\frac{|x-y|^2}{4t}}
u_0(y)\,dy.
\]
It
is shown in \cite[Proposition
4.1]{baoxiang2} that, given
$2\leq p<\infty$,
$1/p+1/p'=1$, $1\leq
q\leq\infty$,
\begin{equation}\label{baox}
\|e^{it\Delta}
u_0\|_{M^{p,q}}\lesssim
(1+|t|)^{-d\left(\frac{1}{2}-\frac{1}{p}\right)}
\|u_0\|_{M^{p',q}}.
\end{equation}
(Similar estimates were
obtained in \cite{benyi}). We
now show that
the condition $p\geq2$ is
necessary in \eqref{baox}, and the decay at
infinity is optimal.
\begin{proposition}{\rm (Sharpness of
\eqref{baox}).}\label{prop3}
Suppose that, for some fixed
$t_0\in\R$, $1\leq
p,q\leq\infty$, $C>0$, the
following estimate holds:
\begin{equation}\label{baox1}
\|e^{it_0\Delta}
u_0\|_{M^{p,q}}\leq C
\|u_0\|_{M^{p',q}},\quad
\forall
u_0\in\mathcal{S}(\R^d).
\end{equation}
Then $p\geq2$.\par
 Assume now
that, for some $\alpha\in\R$,
$C>0$, $M>0$,
$1\leq\gamma,\delta\leq\infty$,
$1\leq p,q\leq\infty$, the
estimate
\begin{equation}\label{baox2}
\|e^{it\Delta}
u_0\|_{M^{p,q}}\leq C
t^\alpha\|u_0\|_{M^{\gamma,\delta}},\quad
\forall
u_0\in\mathcal{S}(\R^d),
\end{equation}
holds for every $t>M$. Then
\begin{equation}\label{baox3}
\alpha\geq
-d\left(\frac{1}{2}-\frac{1}{p}\right).
\end{equation}
\end{proposition}
\begin{proof}
We consider the family of
initial data $u_0(\lambda
x)=e^{-\pi\lambda^2|x|^2}$,
$\lambda>0$. A direct
computation shows that the
corresponding solutions  are
\begin{align}\label{dad}
u(\lambda^2 t,\lambda
x)&=(1+4\pi i t
\lambda^2)^{-d/2}
e^{-\frac{\pi\lambda^2|x|^2}{1+4\pi
i t\lambda^2}}\\
&=\lambda^{-d}\mathcal{G}_{(\lambda^{-2}+4\pi
i t)}(x),\nonumber
\end{align}
where we used the notation in
Lemma \ref{lemm2}. It follows
from \eqref{lem2} that
\begin{equation}\label{ag1}
\|u_0(\lambda\cdot)\|_{M^{p',q}}\asymp
\lambda^{-d}\|\hat{u}_0(\lambda^{-1}\,\cdot)
\|_{W(\Fur
L^{p'},L^q)}=\|\mathcal{\mathcal{G}}_{(\lambda^2)}\|_{W(\Fur
L^{p'},L^q)}\asymp\lambda^{-\frac{d}{p'}},\quad
{\rm as}\ \lambda\to0^+.
\end{equation}
On the other hand, by
\eqref{dad},
\begin{equation}\label{ag2}
\|u(\lambda^2t,\lambda\,\cdot)\|_{M^{p,q}}\asymp
\|\Fur\left(u(\lambda^2t,\lambda\,\cdot)\right)\|_{W(\Fur
L^p,L^q)}\asymp
\lambda^{-d}(a^2+b^2)^{\frac{d}{4}}\|\mathcal{G}_{(a+ib)}\|_{W(\Fur
L^p,L^q)},
\end{equation}
where
\[
a=\frac{\lambda^{-2}}{\lambda^{-4}+(4\pi
t)^2},\quad b=-\frac{4\pi
t}{\lambda^{-4}+(4\pi t)^2}.
\]
Hence, for fixed $t=t_0$,
\eqref{lem2} gives
\begin{equation}\label{ag3}
\|u(\lambda^2
t_0,\lambda\,\cdot)\|_{M^{p,q}}\asymp\lambda^{-\frac{d}{p}},\quad
{\rm as}\ \lambda\to0^+.
\end{equation}
Estimates  \eqref{ag1},
\eqref{ag3} and \eqref{baox1}
yield  $-\frac{d}{p}\geq
-\frac{d}{p'}$, namely
$p\geq2$.\par Choosing $\lambda=1$ in \eqref{ag2} and using
\eqref{lem2}, we obtain
\[
\|u(
t,\cdot)\|_{M^{p,q}}\asymp
t^{-d\left(\frac{1}{2}-\frac{1}{p}\right)},\quad
{\rm as}\ t\to+\infty.
\]
This shows that \eqref{baox3}
is necessary for
\eqref{baox2} to hold.

\end{proof}

\end{document}